\documentclass[12pt]{amsart}
\usepackage{amsfonts, amssymb, amsthm, amsmath, latexsym}

\newtheorem{theorem}{Theorem}[section]
\newtheorem{lemma}[theorem]{Lemma}
\newtheorem{corollary}[theorem]{Corollary}

\newtheorem{definition}[theorem]{Definition}

\theoremstyle{remark}

\newtheorem{remark}[theorem]{Remark}

\newcommand{\func}{\operatorname}

\begin{document}

\title{Hyperbolicity of semigroups and Fourier multipliers}

\author{Yuri Latushkin}
\address{Department of Mathematics\\
University of Missouri\\
Columbia, MO 65211, USA} \email{yuri@math.missouri.edu}
\author{Roman Shvidkoy}
\email{ris1db@mizzou.edu}
\subjclass{Primary:47D06;Secondary:42B15}

\maketitle

\begin{abstract} We present a characterization of hyperbolicity for
strongly continuous semigroups on Banach spaces in terms of Fourier
multiplier properties of the resolvent of the generator. Hyperbolicity
with respect to classical solutions is also considered. Our approach
unifies and simplifies the M. Kaashoek-- S. Verduyn Lunel
theory and multiplier-type
results previously obtained by S. Clark,
M. Hieber, S. Montgomery-Smith,  F. R\"{a}biger, T. Randolph,
and L. Weis.

\end{abstract}

\section{Introduction\label{intro}}

Suppose $X$ is a complex Banach space and $\mathbf{T}=(T_{t})_{t\geq 0}$ is
a strongly continuous semigroup of operators on $X$. Let $A$ denote its
infinitesimal generator.

An autonomous version of a well-known result
     that goes back to O. Perron says the following: a homogeneous
     differential equation $\dot{u}=Au$ admits exponential dichotomy
     on $\mathbb{R}$ if and only if the inhomogeneous equation
     $\dot{u}=Au+f$ has a unique mild solution $u\in F(\mathbb{R};X)$ for
     each $f\in F(\mathbb{R};X)$, see \cite{DK} or \cite{ZL}, and
     \cite{CL}, and the literature therein. Here $F(\mathbb{R};X)$ is
     a space of $X$-valued functions, for instance,
     $F(\mathbb{R};X)=L_{p}(\mathbb{R};X)$, $1\leq p<\infty$. The
     exponential dichotomy for $\dot{u}=Au$ means that the semigroup
     generated by $A$ is
     hyperbolic, that is, condition $\sigma (T_{t})\cap \{
     |z|=1\}=\emptyset$, $t\neq 0$, holds for the spectrum $\sigma (\cdot)$.

     Passing, formally, to the Fourier transforms in the equation
     $\dot{u}=Au+f$ we have that the solution $u$ is given by $u=Mf$,
     where $M:f\mapsto [R(i\cdot ;A)\hat{f}]^{\vee}$, $R(\lambda ;A)$ is
     the resolvent operator, and $\wedge,\vee$ are the Fourier
     transforms. Thus, heuristically, the above-mentioned Perron-type
     theorem could be reformulated to state that the hyperbolicity of
     the semigroup is equivalent to the fact that the function
     $s\mapsto R(is;A)$ is a Fourier multiplier on
     $L_{p}(\mathbb{R};X)$, $1\leq p<\infty$, see, e.g., \cite{Am,Hieber}
     for the definition of Fourier multipliers. One of the objectives
     of the current paper is to systematically study the connections of
     hyperbolicity and $L_{p}$-Fourier multiplier properties of the
     resolvent.

     The use of Fourier multipliers for stability and hyperbolicity for
     strongly continuous semigroups has a fairly long history. To put
     our paper in this context, we briefly review relevant results.
     Probably, the first Fourier multiplier type result was obtained
     in the important paper \cite{KVL} by M. Kaashoek and S.
     Verduyn Lunel. These authors used scalar functions (``matrix
     elements'' of the resolvent) defined by
     $$r_{\rho}(s,x,x^{*})=\langle x^{*}, R(\rho+is;A) x\rangle ,
     \quad \rho\in
     \mathbb{R},\quad s\in \mathbb{R},\quad x\in X,
     \quad x\in X^{*}.$$
     They proved that {\bf T} is hyperbolic if and only if the
     following two conditions holds:
     \begin{enumerate}\item[(i)] $|\langle r_{\rho},\Phi \rangle |\leq K\|
     x\| \|x^{*}\| \| \check{\Phi} \|_{L^{1}}$
     for some $K>0$, $\rho_{0}>0$ and all $\rho$ with $|\rho|<\rho_{0}$ and all
     $\Phi \in \mathcal{S}$, the Schwartz class of scalar functions on
     $\mathbb{R}$;
     \item[(ii)] the C\'esaro integral
     \begin{equation*}\begin{split}
     G_{0}x&=\frac{1}{2\pi}(C,1)\int_{\mathbb{R}}R(is;A)xds\\
     &=\frac{1}{2\pi}\lim_{N\to
        \infty}\frac{1}{N}\int^{N}_{0}\int^{l}_{-l}R(is,A)x ds dl
    \end{split}\end{equation*}
    converges for all $x\in X$.\end{enumerate}

    \noindent Remark, that one of the results of the
    current paper (Theorem~\ref{t5}) shows that condition (ii), in
    fact, follows from (i).

    L. Weis in \cite{WeMnlt} used Fourier multiplier properties of the
    resolvent on Besov spaces to give an alternative proof of the fact
    that the supremum $\omega_\alpha({\bf T})$
    of the growth bounds of ``$\alpha$--smooth'' solutions $T_{t}x$ are
    majorated by the boundedness abscissa $s_{0}(A)$ of the resolvent.
    Remark, that in Section 3 of the present paper we derive a
    formula (Theorem \ref{constants})
    for $\omega_\alpha({\bf T})$ in terms of Fourier multipliers on $L_{p}$.
    Moreover, in Sections 4 and 5 we use Fourier multipliers to study
    an analogue of dichotomy (hyperbolicity) for the smooth solutions.

    A similar formula for $\omega_0({\bf T})$ in terms of the
    resolvent of the generator was obtained in \cite{CLMSR}, see
    also \cite{LM} and formula
    (5.20) in \cite{CL}.
    Formally, Fourier multipliers have not been used in \cite{LM}
    and \cite{CLMSR}. The hyperbolicity of {\bf T} was characterized
    in \cite{LM} and \cite{CLMSR}, see also \cite{CL}, in terms of the
    invertibility of generator $\Gamma$ of the evolution semigroup $\{
    E^{t}\}$ defined on $L_p(\mathbb{R};X)$ as
    $(E^tf)(\tau )=e^{tA}f(\tau -t)$.  However, a simple calculation
    (see Remark~2.2 below) shows that $\Gamma^{-1}=-M$. Thus, formula
    (5.20) in \cite{CL} for the growth bound of {\bf T} is, in fact, a
    Fourier multiplier result that is generalized in
    Theorem \ref{constants} below.

    Via completely different approach based
    on an explicit use of Fourier multipliers, M. Hieber \cite{Hi2}
    gave a characterization of uniform stability for {\bf T} in terms
    of Fourier multiplier properties
    of the resolvent. Also, he proved a formula for $\omega_0({\bf T})$
    that is contained in Theorem \ref{constants} when $\alpha=0$. An
    important ingredient of his proof was the use of well-known
     Datko-van Neerven Theorem saying that
    {\bf T} is uniformly stable if and only if the convolution with
    {\bf T} is a bounded operator on $L_{p}(\mathbb{R};X)$. Since
    the resolvent is the Fourier transform of {\bf T}, the Fourier
    multipliers characterization of uniform stability follows.

    Among other things, this result with a different proof was
    given in \cite{LR}, where Datko-van Neerven Theorem was also used.
    In fact, Theorem \ref{constants} was proved in \cite{LR} for
    $\alpha=0$ or $\alpha=1$. Also, a spectral mapping theorem from \cite{LM} was
    explained in \cite{LR} using Fourier multipliers instead of
    evolution semigroups. In addition, a particular case of Theorem
    \ref{discretealpha} of the current paper (with a different proof) was established
    in \cite{LR}.
    Thus, in the
    present paper we use new technique to ``tie the ends'',
    and give a universal treatment for the results in
    \cite{KVL,CLMSR,Hi2,LR} in a more general context.

    {\bf Acknowledgment.}\quad Yuri Latushkin was supported by the
    Summer Research Fellowship and by the Research Board of the
    University of Missouri. He thanks S. Verduyn Lunel for
    fruitful discussions during his visit to Amsterdam; without these
    discussions this paper would has not been written. Roman Shvidkoy
    was partially supported by the NSF grant DMS-9870027.

\section{\strut Characterization of hyperbolicity\label{char}}

Let us fix some notation:
\begin{itemize}
\item  $\mathbf{T=}(T_{t})_{t\geq 0}$ is a strongly continuous semigroup on
a Banach space $X$ with the generator $A$;
\item $\mathcal{L}(X)$ -- the set of bounded linear operators on $X$;
\item  $R(\lambda ,A)=R(\lambda )$ is the resolvent of $A$;
\item  $\omega_0=\omega _{0}(\mathbf{T)}$ denotes the growth bound of $\mathbf{T}$,
i.e. $\omega_0({\bf T})=\inf \{\omega \in \mathbb{R}:\Vert T_{t}\Vert \leq M_{\omega }e^{\omega
t}\}$;
\item  $s_{0}(A)$ denotes the abscissa of uniform boundedness of the resolvent,
i.e. $s_0(A)=\inf \left\{ s\in \mathbb{R}:\sup \{\Vert R(\lambda )\Vert :\func{Re}%
\lambda >s\}<\infty \right\}$;
\item  $r_{\rho }(s,x,x^{*})=r_{\rho }(s)=\langle x^{*},R(is+\rho )x\rangle$
; $s\in \mathbb{R}$, $x\in X$, $x^{*}\in X^{*}$, $\rho \in \mathbb{R}$;
\item  $\widehat{f}(t)=\int_{\mathbb{R}}f(s)e^{-ist}ds$; $\check{f}(t)=\frac{1}{%
2\pi }\int_{\mathbb{R}}f(s)e^{ist}ds$;
\item  $\mathcal{S}$ stands for the class of Schwartz functions;
\item $\langle r,\Phi \rangle$ denotes the value of a distribution $r$
on $\Phi\in\mathcal{S}$.
\end{itemize}

\begin{definition}
\label{dhyper}We say that the semigroup $\mathbf{T}$\ is \textit{hyperbolic%
} if there is a bounded projection $P$\strut ~on $X$, called splitting,
such that $PT_{t}=T_{t}P$ for all $t>0$ and there exist positive numbers $%
\omega $ and $M$ such that

\begin{enumerate}
\item  $\Vert T_{t}x\Vert \leq Ke^{-\omega t}\Vert x\Vert $, for all $t>0$
and $x\in \func{Im}P$,
\item  $\Vert T_{t}x\Vert \geq Ke^{\omega t}\Vert x\Vert $, for all $t<0$
and $x\in \func{Ker}P$.
\end{enumerate}
The semigroup $\mathbf{T}$ is called uniformly exponentially stable if $P=I$.
\end{definition}

In other words, conditions 1 and 2 say that $(T_{t})_{t\geq 0}$ is uniformly
exponentially stable on $\func{Im}P$ ; all the $T_{t}$'s are invertible on $%
\func{Ker}P$ and the semigroup $(T_{-t})_{t\geq 0}$ is uniformly
exponentially stable there.

\begin{definition}
\label{Greenfn}The function
\begin{equation*}
G(t)=\left\{
\begin{array}{ll}
T_{t}P, & t>0 \\
-T_{t}(I-P), & t<0
\end{array}
\right.
\end{equation*}
is called \textit{the Green's function} corresponding to the hyperbolic
semigroup $\mathbf{T}$.
\end{definition}

Definition \ref{dhyper} allows an equivalent reformulation in terms of
spectral properties of $\mathbf{T}$. Namely, $\mathbf{T}$\ is hyperbolic if
and only if the unit circle $\mathbb{T}$ lies in the resolvent set of $T_{t}$
for one/all $t$ (see \cite[Proposition V.1.15]{Engel}).

Let us recall the following inversion result.
\begin{lemma}
\label{inv}Suppose $\rho >s_{0}(A)$ and $x\in X$, then
\begin{equation*}
F_{t}(x)=\frac{1}{2\pi i}(C,1)\int_{\func{Re}\lambda =\rho }e^{\lambda
t}R(\lambda )xd\lambda \text{, \ \ \ }t\in \mathbb{R}\text{,}
\end{equation*}
where $F_{t}$ is defined as
\begin{equation*}
F_{t}(x)=\begin{cases} T_tx, & t>0\\ \frac{1}{2} x, & t=0\\ 0, &
t<0\end{cases}\text{.}
\end{equation*}
In particular, $\check{r}_{\rho }(t,x,x^{*})=e^{-\rho t}\langle
x^{*},F_{t}(x)\rangle .$
\end{lemma}

The proof can be found in \cite[Theorem 1.3.3]{V-N}. See also Corollary \ref
{inversion}.

Below we establish some algebraic properties of the distributions $\check{r}%
_{\rho }$ for small $|\rho| $ without any additional assumptions on $s_{0}(A)$. The
reader will easily recognize the semigroup properties, in the case $
s_{0}(A)<0$.

In order to be able to threat $r_{\rho }$'s as distributions and to justify
some computations, we assume that the function $s\mapsto \Vert R(is)\Vert $ is bounded
on $\mathbb{R}$, though the proofs below require merely that this function grows not
faster then a power of $|s|$.

\begin{lemma}
\label{l1} If $\tau >0$, then
\begin{equation*}
\check{r}_{0}(t-\tau ,T_{\tau }x,x^{*})=\check{r}_{0}(t,x,x^{*})-\langle
x^{*},T_{t}x\rangle \chi _{[0,\tau ]}(t),\quad t\in \mathbb{R}.
\end{equation*}
\end{lemma}

\begin{proof}
Let us take arbitrary $\Phi \in \mathcal{S}$. Then
\begin{multline*}
\langle \check{r}_{0}(\cdot -\tau ,T_{\tau }x,x^{*}),\Phi \rangle =\langle
\check{r}_{0}(\cdot ,T_{\tau }x,x^{*}),\Phi (\cdot +\tau )\rangle \\
=\langle r_{0}(\cdot ,T_{\tau }x,x^{*}),e^{-i\tau \cdot }\check{\Phi}\rangle
=\int_{-\infty }^{+\infty }\langle x^{*},R(is)e^{-is\tau }T_{\tau }x\rangle
\check{\Phi}(s)ds.
\end{multline*}
Note that
\begin{equation}
e^{-is\tau }R(is)T_{\tau }x=R(is)x-\int_{0}^{\tau }T_{r}x\cdot e^{-isr}dr.
\label{eq1}
\end{equation}
Continuing the line of equalities, we obtain:
\begin{multline*}
\langle \check{r}_{0}(\cdot -\tau ,T_{\tau }x,x^{*}),\Phi \rangle
=\int_{-\infty }^{+\infty }\langle x^{*},R(is)x\rangle \check{\Phi}(s)ds \\
-\int_{-\infty }^{+\infty }\left\langle x^{*},\int_{0}^{\tau
}T_{r}e^{-isr}xdr\right\rangle \check{\Phi}(s)ds=\langle \check{r}_{0},\Phi
\rangle \\
-\int_{0}^{\tau }\langle x^{*},T_{r}x\rangle \Phi (r)dr=\langle \check{r}%
_{0}-\langle x^{*},T_{\cdot }x\rangle \chi _{[0,\tau ]}(\cdot ),\Phi \rangle
.
\end{multline*}
\end{proof}

\begin{lemma}
\label{l2} $\check{r}_{0}(t,T_{\tau }x,x^{*})=\check{r}_{0}(t,x,T_{\tau
}^{*}x^{*})$, $\tau >0$, $t\in \mathbb{R}$.
\end{lemma}

The proof is obvious.

\begin{lemma}
\label{l3} $\check{r}_{0}(t)=e^{\rho t}\check{r}_{\rho }(t)$ for all $\rho $
with $|\rho |<\rho _{0}$ and $t\in \mathbb{R}$.
\end{lemma}

\begin{proof}
We choose $\rho _{0}$ such that $\sup \{\Vert R(is+\rho )\Vert :s\in %
\mathbb{R},|\rho |<\rho _{0}\}$ is finite. Suppose $\Phi \in \mathcal{S}$
has compact support. Then $\check{\Phi}$ is an entire function. Moreover,
\begin{equation}
\lim_{\alpha \to \infty }\check{\Phi}(\alpha +i\beta )=0\text{,}
\label{limit}
\end{equation}
uniformly for all $\beta $ from some finite interval $[a,b]$. It is an
immediate consequence of the following equality:
\begin{equation*}
\check{\Phi}(\alpha +i\beta )=\int_{\mathbb{R}}e^{\beta x}\Phi (x)e^{i\alpha
x}dx=-\frac{1}{i\alpha }\int_{\mathbb{R}}[\beta e^{\beta x}\Phi (x)+e^{\beta
x}\Phi ^{\prime }(x)]e^{i\alpha x}dx.
\end{equation*}
Now using Cauchy's theorem and (\ref{limit}) we get
\begin{equation*}
\begin{split}
\langle \check{r}_{\rho },\Phi \rangle & =\langle r_{\rho },\check{\Phi}%
\rangle =\int_{\mathbb{R}}\langle x^{*},R(is+\rho )x\rangle \check{\Phi}(s)ds
\\
& =\int_{\mathbb{R}}\langle x^{*},R(i(s-i\rho ))x\rangle \check{\Phi}%
(s)ds=\int_{\mathbb{R}-i\rho }\langle x^{*},R(i\lambda )x\rangle \check{\Phi}%
(\lambda +i\rho )d\lambda \\
& =\int_{\mathbb{R}}\langle x^{*},R(is)x\rangle \check{\Phi}(s+i\rho )ds \\
& \quad +2i\lim_{\alpha \to \pm \infty }\int_{0}^{\rho }\langle
x^{*},R(i(\alpha +i\beta ))x\rangle \check{\Phi}(\alpha +i(\beta +\rho
))d\beta \\
& =\int_{\mathbb{R}}\langle x^{*},R(is)x\rangle \check{\Phi}(s+i\rho
)ds=\langle r_{0},\check{\Phi}(\cdot +i\rho )\rangle \\
& =\langle \check{r}_{0},e^{-\rho \cdot }\Phi \rangle ,
\end{split}
\end{equation*}
and the result follows.
\end{proof}

\strut Now we are in a position to prove our main theorem. Let us denote by $M_{\rho }$ the operator
acting by the rule
$$
{M_\rho: f \mapsto [R(i\cdot +\rho )\hat{f}]^{\vee }.}$$
Recall that a function $m\in L_\infty(\mathbb{R};\mathcal{L}(X))$ is
called a {\em Fourier multiplier} on $L_p(\mathbb{R}; X)$ if the
operator $M: f \mapsto [m(\cdot)\hat{f}]^{\vee }$ is a bounded operator
on $L_p(\mathbb{R}; X)$.
Let $L_{1,\infty }(\mathbb{R};X)$ denote the weak-$L_{1}$ space with values
in $X$ (see, e.g., \cite[1.18.6]{Trieb}), that is, the set of all
$X$-valued strongly continuous functions $f$ with the finite norm
\[
\Vert f\Vert _{L_{1,\infty }} :=\sup_{\sigma >0}
\left\{ \sigma \text{mes }
\Big(\{s\in \mathbb{R}:\Vert f(s)\Vert
\ge \sigma \}\Big)<\infty \right\} .\]
Note that $L_{1,\infty }(\mathbb{R};X)\subset L_{1}(\mathbb{R};X)$.

\begin{theorem}
\label{t5} For a strongly continuous semigroup $\mathbf{T}$\ on $X$ the
following conditions are equivalent:
\begin{enumerate}
\item[1)]  $\mathbf{T}$ is hyperbolic;
\item[2)]  $R(i\cdot )$ is a Fourier multiplier on $L_{p}(\mathbb{R},X)$ for
some/all $p$, $1\leq p<\infty $;
\item[3)]  There exists a $\rho _{0}>0$ such that for all $\rho $ with $%
|\rho |<\rho _{0}$, $M_{\rho }$ maps $L_{1}(\mathbb{R},X)$ into $L_{1,\infty
}(\mathbb{R},X)$;
\item[4)]  There exists a $\rho _{0}>0$ such that for all $\rho $ with $%
|\rho |<\rho _{0}$ and all $\Phi \in \mathcal{S}$ we have
$|\langle r_{\rho },\Phi \rangle |\leq K_{\rho }\Vert x\Vert \Vert x^{*}\Vert
\Vert \check{\Phi}\Vert _{1}$.
\end{enumerate}
\noindent Furthermore, if one of these properties holds, then for every $t\in %
\mathbb{R}$ and $x\in X$ the integral
\begin{equation*}
G(t)x=\frac{1}{2\pi }(C,1)\int_{\mathbb{R}}R(is)xe^{ist}ds
\end{equation*}
converges and represents the Green's function of $\mathbf{T}$. Moreover,
$M_{0}f=G*f$ for $f\in L_{1}(\mathbb{R},X)$,
and the splitting projection is given by the formula
\begin{equation}\label{3proj}
P=\frac{1}{2}I+G(0).
\end{equation}
\end{theorem}
\begin{proof}
1)$\Rightarrow $4). This is a part of Theorem 0.2 from \cite{KVL}.

4)$\Rightarrow $2). It follows from 4) that
$\check{r}_{\rho }\in L_{\infty },|\rho |<\rho _{0}$ and
$\Vert \check{r}_{\rho }\Vert _{\infty }\leq K_{\rho }\Vert x\Vert \Vert
x^{*}\Vert$.
By Lemma~\ref{l3}, $\check{r}_{0}(t)=e^{-\rho t}\check{r}_{-\rho }(t)$ a.e.
and $\check{r}_{0}(t)=e^{\rho t}\check{r}_{\rho }(t)$ a.e. for some $\rho >0$%
. So, $|\check{r}_{0}(t)|\leq e^{-\rho |t|}K\Vert x\Vert \Vert x^{*}\Vert $
a.e. for every $x\in X$ and $x^{*}\in X^{*}$, where $K=\max \{K_{\rho
},K_{-\rho }\}$. Now let us fix $p$, $1\leq p<\infty $, and consider a
function $\Phi =\sum_{k=1}^{n}x_{k}\otimes \Phi _{k}$, where $\Phi _{k}\in
\mathcal{S}$ and $\{\Phi _{k}\}$ have disjoint supports. Then $\Vert \Phi
\Vert _{L_{p}}^{p}=\sum_{k=1}^{n}\Vert x_{k}\Vert ^{p}\Vert \Phi _{k}\Vert
_{L_{p}}^{p}$. So, we get the following estimates:
{\allowdisplaybreaks
\begin{eqnarray*}
\Vert M_{0}(\Phi )\Vert _{L_{p}}^{p}& = &\int_{\mathbb{R}}\Vert M_{0}(\Phi
)(t)\Vert ^{p}dt=\frac{1}{2\pi }\int_{\mathbb{R}}\left\| \int_{\mathbb{R}%
}R(is)\hat{\Phi}(s)e^{ist}ds\right\| ^{p}dt \\
& =&\frac{1}{2\pi }\int_{\mathbb{R}}\sup_{\Vert x^{*}\Vert \leq 1}\left|
\sum_{k=1}^{n}\int_{\mathbb{R}}r_{0}(s,x_{k},x^{*})\hat{\Phi}%
_{k}(s)e^{ist}ds\right| ^{p}dt \\
& =&\frac{1}{2\pi }\int_{\mathbb{R}}\sup_{\Vert x^{*}\Vert \leq 1}\left|
\sum_{k=1}^{n}\int_{\mathbb{R}}\check{r}_{0}(\tau ,x_{k},x^{*})\Phi
_{k}(t-\tau )d\tau \right| ^{p}dt \\
& \leq & K^{p}\int_{\mathbb{R}}\left( \sum_{k=1}^{n}\int_{\mathbb{R}}e^{-\rho
|\tau |}\Vert x_{k}\Vert |\Phi _{k}(t-\tau )|d\tau \right) ^{p}dt \\
& = & K^{p}\int_{\mathbb{R}}\left( \int_{\mathbb{R}}e^{-\rho |\tau |}\left(
\sum_{k=1}^{n}\Vert x_{k}\Vert |\Phi _{k}(t-\tau )|\right) d\tau \right)
^{p}dt \\
& \leq & C_{\rho }K^{p}\int_{\mathbb{R}}\int_{\mathbb{R}}e^{-\rho |\tau
|}\left( \sum_{k=1}^{n}\Vert x_{k}\Vert |\Phi _{k}(t-\tau )|\right)
^{p}d\tau dt \\
& = & C_{\rho }K^{p}\int_{\mathbb{R}}e^{-\rho |\tau |}\int_{\mathbb{R}}\left(
\sum_{k=1}^{n}\Vert x_{k}\Vert |\Phi _{k}(t-\tau )|\right) ^{p}dtd\tau  \\
& = & C_{\rho }K^{p}\int_{\mathbb{R}}e^{-\rho |\tau |}\sum_{k=1}^{n}\Vert
x_{k}\Vert ^{p}\Vert \Phi _{k}\Vert _{L_{p}}^{p}d\tau =C_{\rho }^{\prime
}K^{p}\Vert \Phi \Vert _{L_{p}}^{p}.
\end{eqnarray*}}
Since the functions $\Phi $ are dense in $L_{p}(\mathbb{R},X)$, the proof of
4)$\Rightarrow $2) is finished.

2)$\Rightarrow $1). Suppose 2) holds for some $p$, $1\leq p<\infty $. Then,
by the transference principle (see, for example, \cite[Thm. VII.3.8]{StWe}),
$\{R(ik+i\xi )\}_{k\in \mathbb{Z}}$ is a multiplier in $L_{p}(\mathbb{T},X)$
for all $\xi \in \mathbb{R}$, where $\mathbb{T}$ is the unit circle. So,
using results from \cite[Theorem 2.3]{LM} or \cite[Theorem 1]{LR}, we
conclude that $e^{2\pi i\xi }\in \rho (T_{2\pi })$ for all $\xi \in %
\mathbb{R}$. Thus, $\mathbb{T}\subset \rho (T_{2\pi })$ and hence $\mathbf{T}
$\ is hyperbolic.

This completes the proof of 1)$\Leftrightarrow $2)$\Leftrightarrow $4).

2)$\Rightarrow $3). It is easy to see using the resolvent identity, that
there exists a $\rho _{0}>0$ such that $R(i\cdot +\rho )$ is a $L_{1}(%
\mathbb{R},X)$-multiplier for all $\rho $ such that $|\rho |<\rho _{0}$.

3)$\Rightarrow $4). Without loss of generality, assume $\rho =0$. Denote
\begin{equation*}
\mu=\sup_{0\leq \tau \leq 1}\Vert T_{\tau }\Vert
\end{equation*}
and fix $x\in X$, $x^{*}\in X^{*}$, $\Vert x\Vert =\Vert x^{*}\Vert =1$. Let
us take a function $\Phi \in \mathcal{S}$. By condition 3) we have
\begin{equation*}
\Vert M_{0}(\check{\Phi}\otimes x)\Vert _{1,\infty }\leq K\Vert \hat{\Phi}%
\Vert _{1}.
\end{equation*}
So, $\func{mes}\{\tau :\Vert M_{0}(\hat{\Phi}\otimes x)(\tau )\Vert >2K\Vert
\hat{\Phi}\Vert _{1}\}\leq \frac{1}{2}$. This implies that there is a $\tau $%
, $-1<\tau <0$, such that
\begin{equation*}
\Vert M_{0}(\check{\Phi}\otimes x)(\tau )\Vert \leq 2K\Vert \hat{\Phi}\Vert
_{1}.
\end{equation*}
Let us apply the functional $T_{-\tau }^{*}x^{*}$ to the left-hand side of
the inequality. Then we have:
\begin{equation*}
|\frac{1}{2\pi }<T_{-\tau }^{*}x^{*},[R\cdot \Phi \otimes x]^{\vee }(\tau
)>|\leq 2\mu KC\Vert \hat{\Phi}\Vert _{1}.
\end{equation*}
By Lemma~\ref{l1} and \ref{l2} the expression under the absolute value sign
is equal to
\begin{multline*}
\check{r}_{0}(\cdot ,x,T_{-\tau }^{*}x^{*})*\check{\Phi}(\tau )=\langle
\check{r}_{0}(\cdot ,T_{-\tau }x,x^{*}),\check{\Phi}(\tau -\cdot )\rangle  \\
=\langle \check{r}_{0}(\cdot +\tau ,T_{-\tau }x,x^{*}),\hat{\Phi}\rangle
=\langle \check{r}_{0}(\cdot ,x,x^{*}),\hat{\Phi}\rangle -\int_{0}^{-\tau
}\langle x^{*},T_{t}x\rangle \hat{\Phi}(t)dt.
\end{multline*}
By the triangle inequality, we have
\begin{equation*}
|\langle r_{0}(\cdot ,x,x^{*}),\Phi \rangle |\leq 2K\mu\Vert \hat{\Phi}\Vert
_{1}+\mu\Vert \hat{\Phi}\Vert _{1}\leq 3K\mu\Vert \check{\Phi}\Vert _{1},
\end{equation*}
which is what we wanted.

Now we turn to the second part of the theorem. First, we prove an auxiliary
Fej\'er-type lemma (probably, well-known).

\begin{lemma}
\label{l6} If $f\in L_{1}(\mathbb{R},X)$, then the integral
\begin{equation*}
\frac{1}{2\pi }(C,1)\int_{\mathbb{R}}\hat{f}(s)e^{ist}ds
\end{equation*}
converges to $f(t)$ a.e. Moreover,
\begin{equation*}
f=\frac{1}{2\pi }L_{1}-\lim_{N\to \infty }\frac{1}{N}\int_{0}^{N}\int_{-\ell
}^{\ell }\hat{f}(s)e^{is\cdot }ds.
\end{equation*}
\end{lemma}

\begin{proof}
\begin{equation*}
\begin{split}
\frac{1}{2\pi}\frac{1}{N} \int^N_0\int^\ell_{-\ell}\hat{f}(s)e^{ist}ds&=%
\frac{1}{2\pi}\int^N_{-N}\hat{f}(s)e^{ist}\left( 1-\frac{|s|}{N}\right) ds \\
&= \frac{1}{2\pi}\int^N_{-N}\int^{+\infty}_{-\infty}f(r)e^{isr}dr\cdot
e^{ist}\left( 1-\frac{|s|}{N}\right) ds \\
&= \int^{+\infty}_{-\infty}f(r)\frac{1}{2\pi}\int^N_{-N}e^{is(t-r)}\left( 1-%
\frac{|s|}{N}\right) dsdr.
\end{split}
\end{equation*}
The inner integral is equal to $K_N(t-r)=\frac{1}{\pi N(t-r)^2}[1-\cos N(t-r)%
]$. One can easily check that $K_N$ is a positive kernel in $L_1$, that is, $%
(K_N*f)(\cdot)$ tends to $f(\cdot)$ a.e. and in $L_1$ as $N\to\infty$.
\end{proof}

Suppose $f\in L_{1}(\mathbb{R})$. Then by 2) we have that $M_{0}(f\otimes x)\in
L_{1}(\mathbb{R},X)$. By Lemma~\ref{l6}, there is a $\tau \in (-1,0)$ such
that
\begin{equation*}
M_{0}(f\otimes x)(\tau )=\frac{1}{2\pi }(C,1)\int_{\mathbb{R}}R(is)x\hat{f}%
(s)e^{is\tau }ds.
\end{equation*}
Let us apply the operator $T_{-\tau }$. Then using \eqref{eq1} we obtain:
\begin{equation*}
\begin{split}
T_{-\tau }([R\hat{f}\otimes x]^{\vee }(\tau ))& =\frac{1}{2\pi }(C,1)\int_{%
\mathbb{R}}R(is)T_{-\tau }xe^{is\tau }\hat{f}(s)ds \\
& =\frac{1}{2\pi }(C,1)\int_{\mathbb{R}}\Big[ R(is)x\hat{f}(s) \\
& \qquad -\int_{0}^{-\tau }T_{r}xe^{-isr}dr\cdot \hat{f}(s)\Big]ds.
\end{split}
\end{equation*}
Since
\begin{equation*}
f(-r)=L_{1}-\lim_{N\to \infty }\frac{1}{2\pi }\frac{1}{N}\int_{0}^{N}\int_{-%
\ell }^{\ell }\hat{f}(s)e^{-isr}dsd\ell
\end{equation*}
and $V\varphi =\int_{0}^{-\tau }T_{r}x\cdot \varphi (r)dr$ is a bounded
linear operator from $L_{1}(\mathbb{R})$ to $X$, we conclude that the $(C,1)$%
-integral of the second summand converges and equals $\int_{0}^{-\tau
}T_{r}x\cdot f(-r)dr$. This means, in particular, that $\frac{1}{2\pi }%
(C,1)\int_{\mathbb{R}}R(is)x\hat{f}(s)ds$ converges. Let us denote it by $%
G(0,f)$. Also let
\begin{equation*}
G(t,f)=G(0,f(\cdot -t))=\frac{1}{2\pi }(C,1)\int_{\mathbb{R}}R(is)x\hat{f}%
(s)e^{ist}ds
\end{equation*}
for $f\in L_{1}$, $t\in \mathbb{R}$, $x\in X$. Now we introduce the
following operators:
\begin{align}
S_{N}^{t}(f,x)& =\frac{1}{2\pi }\frac{1}{N}\int_{0}^{N}\int_{-\ell }^{\ell
}R(is)x\hat{f}(s)e^{ist}dsd\ell ; \nonumber\\
I_{N}^{t}(x)& =\frac{1}{2\pi }\frac{1}{N}\int_{0}^{N}\int_{-\ell }^{\ell
}R(is)xe^{ist}dsd\ell .\label{defSNt}
\end{align}
It is easy to see that $\Vert S_{N}^{t}(f,x)\Vert \leq C_{N}\Vert f\Vert
_{L_{1}}\Vert x\Vert $. On the other hand, we have just proved that $%
G(t,f)x=\lim_{N\to \infty }S_{N}^{t}(f,x)$ exists for all $f\in L_{1}$, $%
x\in X$. So, by the boundedness principle for bilinear operators, $\Vert
S_{N}^{t}\Vert \leq C$, where $C$ does not depend on $N$ and $t$.

Let $f_{\epsilon }$, $\epsilon >0$, be a kernel in $L_{1}(%
\mathbb{R})$, that is, $f_{\epsilon }*\Phi \to \Phi $ as $\epsilon \to 0$
for each $\Phi \in L_{1}(\mathbb{R})$. Then $I_{N}^{t}(x)=\lim_{\epsilon \to
0}S_{N}^{t}(f_{\epsilon },x)$ and hence, $\Vert I_{N}^{t}\Vert \leq C$.

Let us show that $G(t)x=\lim_{N\to \infty }I_{N}^{t}(x)$ exists for all $%
x\in D(A^{2})$. This will be enough to prove that $G(t)x$ exists for all $%
x\in X$. Fix $x\in D(A^{2})$ and notice that
\begin{equation*}
\begin{split}
& I_N^t(x)=\frac{1}{2\pi}\frac{1}{N}\int\limits_0^N
\int\limits_{-\ell}^{\ell}R(is)xe^{ist}dsd\ell\\
&=\frac{1}{2\pi}\frac{1}{N}\int\limits_0^1
\int\limits_{-\ell}^{\ell}R(is)xe^{ist}dsd\ell
+\frac{1}{2\pi}\frac{1}{N}\int\limits_1^N
\int\limits_{|s|\le 1}R(is)xe^{ist}dsd\ell\\
&\qquad\qquad +\frac{1}{2\pi}\frac{1}{N}\int\limits_1^N
\int\limits_{1\le|s|\le\ell}R(is)xe^{ist}dsd\ell\\
&=\frac{1}{2\pi}\frac{1}{N}\int\limits_0^1
\int\limits_{-\ell}^{\ell}R(is)xe^{ist}dsd\ell
+\frac{1}{2\pi}\frac{N-1}{N}\int\limits_{|s|\le 1}R(is)xe^{ist}ds\\
&\qquad\qquad
+\frac{1}{2\pi }\frac{1}{N}\int_{1}^{N}\int_{1\leq |s|\leq \ell
}\left[ -\frac{R(is)A^{2}x}{s^{2}}+\frac{x}{is}-\frac{Ax}{s^{2}}\right]
e^{ist}dsd\ell .
\end{split}
\end{equation*}
So,
\begin{eqnarray}\label{int}
\lim_{N\to \infty }I_{N}^{t}(x) &=&\frac{1}{2\pi }\int_{|s|\leq 1}R(is)xe^{ist}ds  \notag \\
&-&\frac{1}{2\pi }\int_{|s|\geq 1}\left[ \frac{Ax}{s^{2}}+\frac{R(is)A^{2}x}{s^{2}}\right]
e^{ist}ds \\
&+&\frac{x}{2\pi i}\int_{t}^{+\infty }\frac{sin(s)}{s}ds\cdot \chi _{\mathbb{R} \backslash \{0\}}(t). \notag
\end{eqnarray}
%

Finally, it is only left to verify that $G(t)$ is indeed the Green's
function. Let us prove the first equality in Definition \ref{Greenfn}, the
second one being analogous. We have
\begin{align*}
T_{\tau }Px& =\frac{1}{2}T_{\tau }x+\frac{1}{2\pi }(C,1)\int_{\mathbb{R}%
}R(is)xT_{\tau }xds \\
& =\frac{1}{2}T_{\tau }x+\frac{1}{2\pi }(C,1)\int_{\mathbb{R}}R(is)xe^{ist}ds
\\
& \qquad -\frac{1}{2\pi }(C,1)\int_{\mathbb{R}}\int_{0}^{\tau
}T_{r}xe^{-isr}dr\cdot e^{is\tau }ds \\
& =G(\tau )x+\frac{1}{2}T_{\tau }x-\frac{1}{2}T_{\tau }x=G(\tau )x,
\end{align*}
where we use the ordinary Fej\'{e}r's theorem.

It follows from the above that $G(\tau )$ is an exponentially decaying
function. So, $f\mapsto G*f$ is a bounded operator on $L_{1}$. On the other
hand, $M_{0}\Phi =G*\Phi $ for all $\Phi \in \mathcal{S}$. Hence, $%
M_{0}f=G*f $ for all $f\in L_{1}$.

The proof of \eqref{3proj} can be found in \cite{KVL}.
\end{proof}

\begin{corollary}[\cite{Gear,Herbst}]
Suppose $X$ is a Hilbert space. Then the semigroup $\mathbf{T}$ is
hyperbolic if and only if the resolvent $R(\lambda ,A)$ is bounded in some
strip containing the imaginary axes.
\end{corollary}

\begin{remark}
\label{r7} Condition 3) can be considerably weakened in the following way.
Suppose $F$ is a space of functions on $\mathbb{R}$ with the following
property: for any $f\in F$ there is a $t\in [-1,0]$ such that $|f(t)|\leq
c\Vert f\Vert _{F}$. Many quasi-normed
function spaces have this property, for example, $L_{q,r}(\mathbb{R})$, $H_{p}(\mathbb{R})$, $C_0(\mathbb{%
R)}$, or any function lattice with $\Vert \chi _{[-1,0]}\Vert \neq 0$. Denote
by $F(X)$ the space of all strongly measurable functions $f$ with values in $%
X$ such that $\Vert f(\cdot )\Vert \in F$. Our proof shows that it is enough
to require that $M_{\rho }$ maps $L_{1}(\mathbb{R},X)$ into $F(X)$ (see also
the proof of Theorem \ref{alpha}).
\end{remark}

\begin{remark}
\label{r8} Recall that the generator $\Gamma $ of the evolution semigroup $%
(E^{t})_{t\ge 0}$, defined on $L_{p}(\mathbb{R},X)$ by $%
(E^{t}f)(s)=T_{t}f(s-t)$, is the closure of the operator $-d/dt+A$ on the
domain $D(-d/dt)\cap D(A)$. It is known, see \cite[Thm.2.39]{CL}, that $%
\mathbf{T}$ is hyperbolic if and only if the operator $\Gamma $ is
invertible on one/all $L_{p}(\mathbb{R},X)$, $1\le p<\infty $. This result
immediately implies that conditions 1) and 2) in Theorem \ref{t5} are
equivalent. Indeed, if $x\in D(A)$ and $\Phi \in \mathcal{S}$,
then $\Gamma (\Phi \otimes x)=-\Phi ^{\prime }\otimes x+\Phi
\otimes Ax$. Using elementary properties of Fourier transform, we have that
\begin{equation*}
M_{0}\Gamma (\Phi \otimes x)=-\Phi \otimes x,\quad x\in D(A)
\end{equation*}
and
\begin{equation*}
\Gamma M_{0}(\Phi \otimes x)=-\Phi \otimes x,\quad x\in X\text{,}
\end{equation*}
and the result follows.
\end{remark}

It is worth noting that in the special case $s_{0}(A)<0$, by Lemma \ref{inv}%
, the splitting projection turns into the identity and our theorem gives the
characterization of uniform exponential stability observed in \cite{LR}.

There is a Mikhlin-type sufficient condition due to M. Hieber \cite{Hieber}
for an operator-valued symbol to be $L_{1}$-multiplier. Applied to the
resolvent it yields the following: if there exists a $%
\delta >\frac{3}{4}$ such that $\sup \{|s|^{\delta }\Vert R(is)\Vert
\}<\infty $, then $R(i\cdot )$ is a multiplier.

Yet another condition for operator-valued symbol to be a multiplier is
recently developed in \cite{WeisL}. It works if $X$ is a UMD-space and says
that if the families $\{R(is)\}_{s\in \mathbb{R}}$ and $\{sR^{2}(is)\}_{s\in
\mathbb{R}}$ are R-bounded, then $R(i\cdot )$ is a multiplier.

\section{Extension to the case $\alpha >0\label{ext}$}

It turns out that many arguments from Section \ref{char} work in a more
general situation, when the resolvent multiplier is restricted to $L_{p}(%
\mathbb{R},X_{\alpha })$, where $X_{\alpha }$ is the domain of the fractional
power $(A-\omega )^{\alpha }$, endowed with the norm $\Vert x\Vert _{\alpha
}=\Vert (A-\omega )^{\alpha }x\Vert $.
In this section we show that $R(i\cdot +\rho )$ is a multiplier from $L_{p}(\mathbb{R}%
,X_{\alpha })$ to $L_{p}(\mathbb{R},X)$ for small values of $\rho $ if and only
if the following modified Kaashoek - Verduyn Lunel inequality holds: $%
|\langle r_{\rho },\Phi \rangle |\leq K\Vert x\Vert _{\alpha }\Vert
x^{*}\Vert \Vert \hat{\Phi}\Vert _{L_{1}}$. Also in this case $G(t)x$ exists
for all $x\in X_{\alpha }$ and is exponentially decaying as $|t|\rightarrow
\infty $. As a by-product of this results we obtain the following
relationship between the fractional growth bound $\omega _{\alpha }(\mathbf{T%
})$ and its spectral analogue $s_{\alpha }(A)$ (see \eqref{dsa} for the
definitions): $\omega _{\alpha }(\mathbf{T})$ is the infimum of all
$\omega >s_{\alpha }(A)$ such that $R(i\cdot +\omega )$ is a
multiplier from $L_{p}(\mathbb{R},X_{\alpha })$ to $L_{p}(\mathbb{R},X)$. In the
particular case, when $X$ is a Hilbert space, the latter condition will be
shown to hold for all $\omega >s_{\alpha }(A)$. So, $s_{\alpha }(A)=$ $%
\omega _{\alpha }$, which gives a different proof
 of G. Weiss's \cite{Weiss} result for
 arbitrary $\alpha \geq 0$, also obtained by L. Weis and V. Wrobel in \cite
{WW}.
The main result in this section is an
extension of Theorem \ref{t5} to the case of arbitrary $\alpha >0$. To be
more precise, we treat only conditions 2)-4), as hyperbolicity is ambiguous
in this situation and therefore it is postponed to the next section.

One can notice that most of the proof of Theorem \ref{t5} work for all $%
\alpha >0$ if one replaces all $X$-norms by $X_{\alpha }$-norms. However,
the ``some/all'' part of condition 2), being an easy consequence of
results in \cite{LM} and the spectral characterization of
hyperbolicity in case $\alpha =0$, requires some additional duality argument.

Before we state our main theorem, let us recall the notion of fractional
power of $A$.
Suppose $\omega >\max \{\omega _{0}+3,3\}$. Denote $A-\omega $ by $A_{\omega
}$. Let $\gamma $ be the path consisting of two rays $\Gamma
_{1}=\{-1+te^{i\theta }:t\in [0,+\infty )\}$ and $\Gamma
_{2}=\{-1-te^{i\theta }:t\in [0,+\infty )\}$ going upwards. We assume that $%
\theta $, $\theta <\frac{\pi }{6}$, is small enough to ensure the inequality
$\Vert R(\mu +\omega )\Vert \leq C\frac{1}{1+|\mu |}$ in the sector
generated by $\gamma $. For any $\alpha >0$ we define $A_{\omega }^{\alpha }$
as the inverse to the operator $A_{\omega }^{-\alpha }$ acting on $X$ by the rule
\begin{equation*}
A_{\omega }^{-\alpha }(x)=\frac{1}{2\pi i}\int_{\gamma }\mu ^{-\alpha }R(\mu
+\omega )xd\mu.
\end{equation*}
 Let us denote by $X_{\alpha }$ the domain of $A_{\omega }^{\alpha }$
endowed with the norm $\Vert x\Vert _{\alpha }=\Vert (A-\omega )^{\alpha
}x\Vert $. Then $X_{\alpha }$ is a Banach space and it does not depend on
the particular choice of $\omega $, $\omega >\omega _{0}$,
see \cite{Engel} for more information concerning fractional powers.

\begin{theorem}
\label{alpha} Assume that there exists a $\rho _{0}>0$ such that
\begin{equation}\label{assgrR}
\sup \left\{ \frac{\Vert R(\lambda )\Vert }{1+|\lambda |^{\alpha }}:|\func{%
Re}\lambda |<\rho _{0}\right\} <\infty.
\end{equation}
 Then the following conditions
are equivalent:
\begin{itemize}
\item[1)]  $R(i\cdot +\rho )$ is a multiplier from $L_{p}(\mathbb{R},X_{\alpha
})$ to $L_{p}(\mathbb{R},X)$, for some/all $p$, $1\leq p<\infty $, and all $%
\rho $, $|\rho |<\rho _{0}$;
\item[2)]  $R(i\cdot +\rho )$ is a multiplier from $L_{p}(\mathbb{R},X_{\alpha
})$ to $E(X)$, for some $p$, $1\leq p<\infty $, and all $\rho $, $|\rho
|<\rho _{0}$, where $E$ is a rearrangement invariant quasi-Banach lattice;
\item[3)]  $R(i\cdot +\rho )$ is a multiplier from $L_{1}(\mathbb{R},X_{\alpha
})$ to $F(X)$ for all $\rho $, $|\rho |<\rho _{0}$, where $F$ is some
rearrangement invariant quasi-Banach lattice;
\item[4)]  $|\langle r_{\rho },\Phi \rangle |\leq K\Vert x\Vert _{\alpha
}\Vert x^{*}\Vert \Vert \hat{\Phi}\Vert _{L_{1}}$ for all $\Phi \in \mathcal{%
S}$, $x\in X_{\alpha }$, $x^{*}\in X^{*}$ and $|\rho |<\rho _{0}$.
\end{itemize}
If one of these conditions holds, then the integral
\begin{equation*}
G(t)x=\frac{1}{2\pi }(C,1)\int_{\mathbb{R}}R(is)xe^{ist}ds
\end{equation*}
converges for all $x\in X_{\alpha }$, and $\Vert G(t)x\Vert \leq K\Vert
x\Vert _{\alpha }e^{-\rho |t|}$ for all $\rho $, $0<\rho <\rho _{0}$.
Moreover, $M_{0}f=G*f$ for $f\in L_1(\mathbb{R},X)$.
\end{theorem}

\begin{proof}
1)$\Rightarrow $2) is evident.

2)$\Rightarrow $3). Assume for simplicity that $\rho =0$. First we claim
that $M_{0}$ maps $L_{p}(\mathbb{R}, X_{\alpha })$ into $L_{\infty }(\mathbb{R}, X)$. To prove this,
let us take an arbitrary function $\Phi $ of the form $\sum_{i=1}^{n}\Phi
_{i}x_{i}$ , where $x_{i}\in X_{\alpha }$ and $\Phi _{i}\in \mathcal{S}$.
Then, by condition 2),
\begin{equation*}
\Vert M_{0}(\Phi )\Vert _{E(X)}\leq K\Vert \Phi \Vert _{L_{p}(\mathbb{R}, X_{\alpha })}%
\text{.}
\end{equation*}
It implies that for every $n\in \mathbb{Z}$ there exists a $t\in [n,n+1]$ such
that
\begin{equation*}
\Vert [R\hat{\Phi}]^{\vee }(t)\Vert _{X}\leq \frac{2K}{\varphi (1)}\Vert
\Phi \Vert _{L_{p}(\mathbb{R},X_{\alpha })}\text{,}
\end{equation*}
where $\varphi $ is the characteristic function of $E$. For any fixed $\tau
\in [0,2]$, let us apply the operator $T_{\tau }$ to the right-hand side of
this inequality. Then we get
\begin{equation*}
\Vert \frac{1}{2\pi }\int_{\mathbb{R}}T_{\tau }R(is)\hat{\Phi}(s)e^{ist}ds\Vert
\leq C\Vert \Phi \Vert _{L_{p}(\mathbb{R}, X_{\alpha })}\text{.}
\end{equation*}
Now using equality (\ref{eq1}) we obtain the following
\begin{eqnarray*}
\frac{1}{2\pi }\int_{\mathbb{R}}T_{\tau }R(is)\hat{\Phi}(s)e^{ist}ds &=&\frac{1%
}{2\pi }\int_{\mathbb{R}}R(is)\hat{\Phi}(s)e^{is(t+\tau )}ds \\
&&-\frac{1}{2\pi }\int_{\mathbb{R}}\int_{0}^{\tau }T_{r}\hat{\Phi}(s)e^{is(\tau
-r)}drds \\
&=&M_{0}(\Phi )(t+\tau )-\frac{1}{2\pi }\int_{0}^{\tau }T_{r}\Phi (\tau -r)dr%
\text{.}
\end{eqnarray*}
Thus, $\left\| M_{0}(\Phi )(t+\tau )\right\| _{X}\leq \tilde{C}\Vert \Phi
\Vert _{L_{p}(\mathbb{R}, X_{\alpha })}$. By the choice of $\tau $ and $t$ we have the
same inequality on the whole real line. Since $\tau $ was chosen arbitrary,
the claim is proved.

Let us observe that the boundedness of $M_{0}$ is equivalent to the fact
that $R(i\cdot )A_{\omega }^{-\alpha }$ is an $L_{p}(\mathbb{R}, X)-L_{\infty }(\mathbb{R}, X)$
multiplier.

Denote by $X^{\odot }$ the sun dual to $X$ on which the dual semigroup is
strongly continuous (see \cite{Engel}). One can easily check, by duality,
that for a test function $\Phi =\sum_{i=1}^{n}\Phi _{i}x_{i}^{\odot }$ one
has
\begin{equation*}
\left\| \lbrack R^{\odot }(A^{\odot })^{-\alpha }\hat{\Phi}]^{\vee }\right\|
_{L_{q}(\mathbb{R}, X^{\odot })}\leq \tilde{C}\Vert \Phi \Vert _{L_{1}(\mathbb{R}, X^{\odot })}\text{%
,}
\end{equation*}
where $\frac{1}{p}+\frac{1}{q}=1$ and $A^{\odot }$ is the generator of the
sun dual semigroup. In other words, $M_{0}^{\odot }$ maps $L_{1}(\mathbb{R}, X_{\alpha
}^{\odot })$ into $L_{q}(\mathbb{R}, X^{\odot })$.

By what we just proved, $M_{0}^{\odot }$ is bounded from $L_{1}(\mathbb{R}, X_{\alpha
}^{\odot })$ to $L_{\infty }(\mathbb{R}, X^{\odot })$, and again by duality, $M_{0}$
maps $L_{1}(\mathbb{R}, X_{\alpha })$ into $L_{\infty }(\mathbb{R}, X)$, which proves condition 3)
with $F=L_{\infty }$.

The proofs of all other implications are completely analogous to those of
Theorem \ref{t5}.

Let us now turn to the second part of our theorem. Although its proof is
also essentially the same, some comments will be in order. By Lemma
\ref{growth},
proved below, assumption \eqref{assgrR} is
equivalent to $\sup \left\{ \Vert R(\lambda )\Vert _{X_{\alpha }\rightarrow
X}:|\func{Re}\lambda |<\rho _{0}\right\} <\infty $. So, the operators $%
S_{N}^{t}$, introduced in \eqref{defSNt}, are bounded from $X_{\alpha }\times L_{1}(\mathbb{R})$ to $X$.
Uniform boundedness follows from the fact that $\underset{N\rightarrow
\infty }{\lim }S_{N}^{t}(x,f)$ exists for all $x\in X_{\alpha }$ and $f\in
L_{1}(\mathbb{R})$ by Lemma \ref{l6}. Consequently, $\Vert I_{N}^{t}\Vert
_{X_{\alpha }\rightarrow X}\leq C$. Formula (\ref{int}) still makes sense
for all $x\in X_{\alpha +2}$, because then $A^{2}x\in X_{\alpha }$ and all
the integrals converge absolutely. So, $G(t)x$ exists for all $x\in
X_{\alpha }$, and it is continuous in $t$, $t\neq 0$.

Since $\langle x^{*},G(t)x\rangle =\check{r}_{0}(t,x,x^{*})$, by condition
3) and Lemma \ref{l3} we have that $|\langle x^{*},G(t)x\rangle |\leq Ke^{-\rho
|t|}\Vert x\Vert _{\alpha }\Vert x^{*}\Vert $ almost everywhere and hence,
by the continuity of $G(t)x$, for all $t\in \mathbb{R}$. Thus, $\left\|
G(t)x\right\| _{X}\leq Ke^{-\rho |t|}\Vert x\Vert _{\alpha }$ and the proof
is finished.
\end{proof}

\begin{lemma}
\label{growth} Let $S=\{\lambda \in \mathbb{C}:a<\func{Re}\lambda <b\}$, $a,b\in\mathbb{R}$, be a
subset of $\rho (A)$, where $a\in \mathbb{R}$, $b\in \mathbb{R}$. Then conditions
$$\sup
\left\{ \frac{\Vert R(\lambda )\Vert }{1+|\lambda |^{\alpha }}:\lambda \in
S\right\} <\infty \quad\text{and}\quad \sup \left\{ \Vert R(\lambda )A_{\omega
}^{-\alpha }\Vert :\lambda \in S\right\} <\infty$$
are equivalent.
\end{lemma}

\begin{proof}
Since $b$ is finite, there are constants $c>0$ and $\varphi _{0}$, $%
0<\varphi _{0}<\pi $ such that $|\mu -e^{i\varphi }|>|\mu |+c$ for all $\mu
\in \gamma $ and $\varphi _{0}<|\varphi |<\pi -\varphi _{0}$. Pick
$N>1$ large enough to satisfy $\frac{\omega }{N}<\frac{c}{2}$ and such that
whenever $\lambda \in S$ and $|\lambda |>N$ ,\ then $\ \varphi _{0}<|\arg
\lambda |<\pi -\varphi _{0}$ and $\lambda $ does not belong to the sector
bounded by the contour $|\lambda |\gamma $. For all such $\lambda $ we
have
\begin{equation}
|\mu +\frac{\omega }{|\lambda |}-e^{i\arg \lambda }|>|\mu |+\frac{c}{2}
\label{mu}.
\end{equation}
Let us consider the following integral:
\begin{equation*}
I_{\lambda }=\int_{\gamma }\frac{\mu ^{-\alpha }}{\mu +\omega -\lambda }d\mu
\text{, \ \ }\lambda \in S\text{, \ \ }|\lambda |>N.
\end{equation*}
By the choice of $N$, the integrand does not have singular points between $%
\gamma $ and $|\lambda |\gamma $. By the Cauchy Theorem, we have
\begin{equation*}
I_{\lambda }=\int_{|\lambda |\gamma }\frac{\mu ^{-\alpha }}{\mu +\omega
-\lambda }d\mu =\frac{1}{|\lambda |^{\alpha }}\int_{\gamma }\frac{\mu
^{-\alpha }}{\mu +\frac{\omega }{|\lambda |}-e^{i\arg \lambda }}d\mu \text{.}
\end{equation*}
Inequality (\ref{mu}) implies that the absolute value of the last integral
is bounded from above by a constant that does not depend on $\lambda $, whenever
$\lambda \in S$, \ \ $|\lambda |>N$. The analogous estimate from below
follows from geometric considerations. Thus,
\begin{equation}
\frac{d_{1}}{|\lambda |^{\alpha }}\leq |I_{\lambda }|\leq \frac{d_{2}}{%
|\lambda |^{\alpha }}\text{,}  \label{bounds}
\end{equation}
for some positive $d_{1}$ and $d_{2}$.

Suppose $x\in X$. Then
\begin{eqnarray*}
R(\lambda )A_{\omega }^{-\alpha }x &=&\frac{1}{2\pi i}\int_{\gamma }\mu
^{-\alpha }R(\lambda )R(\mu +\omega )xd\mu \\
&=&\frac{1}{2\pi i}I_{\lambda }R(\lambda )x-\frac{1}{2\pi i}\int_{\gamma }%
\frac{\mu ^{-\alpha }}{\mu +\omega -\lambda }R(\mu +\omega )xd\mu \text{.}
\end{eqnarray*}
Let us notice that $|\mu +\omega -\lambda |\geq K(|\mu |+1)$ for some $K>0$
and all $\mu \in \gamma $, $\lambda \in S$, $|\lambda |>N$, whereas $\Vert
R(\mu +\omega )\Vert \leq C\frac{1}{1+|\mu |}$. Consequently,
\begin{equation*}
\left\| \int_{\gamma }\frac{\mu ^{-\alpha }}{\mu +\omega -\lambda }R(\mu
+\omega )xd\mu \right\| \leq K_{\alpha }\Vert x\Vert \text{.}
\end{equation*}
In combination with \eqref{bounds} this gives the following estimates:
\begin{eqnarray*}
\left\| R(\lambda )A_{\omega }^{-\alpha }x\right\| &\geq &d_{1}\frac{\left\|
R(\lambda )x\right\| }{|\lambda |}-K_{\alpha }\Vert x\Vert, \\
\left\| R(\lambda )A_{\omega }^{-\alpha }x\right\| &\leq &d_{2}\frac{\left\|
R(\lambda )x\right\| }{|\lambda |}+K_{\alpha }\Vert x\Vert \text{,}
\end{eqnarray*}
for all $\lambda \in S$, $|\lambda |>N$ and $x\in X$, which proves the lemma.
\end{proof}

\begin{remark}
In view of Lemma \ref{growth}, assumption \eqref{assgrR} in Theorem \ref{alpha},
in fact, follows from condition 1) or 3).
\end{remark}

\begin{remark}
\label{greenrem}Just as in the proof of Theorem \ref{t5} one can show the
following identities:
\begin{eqnarray}
G(t)=T_{t}P\text{, \ \ }t>0  \label{Green} \\
G(t)T_{-t}=-(I-P)\text{, \ \ \ }t<0  \notag
\end{eqnarray}
on $X_{\alpha }$, where $P$ is defined as $\frac{1}{2}I+G(0)$.
\end{remark}

Let us now recall the definition of the \textit{fractional growth bound} $\omega
_{\alpha }(\mathbf{T})$ and its spectral counterpart $s_{\alpha }(A)$:
\quad $\omega _{\alpha }(\mathbf{T})$ is the infimum of all $\omega \in
\mathbb{R}$ such that $\Vert T_{t}x\Vert \leq M_{\omega }e^{\omega t}\Vert
x\Vert _{\alpha }$, for some $M_{\omega }>0$ and all $x\in X_{\alpha }$ and $%
t\geq 0$, and
\begin{equation}\label{dsa}
s_{\alpha }(A)=\inf \left\{ s:\sup \left\{ \frac{\Vert R(\lambda
)\Vert }{1+|\func{Im}\lambda |^{\alpha }}:\func{Re}\lambda >s\right\}
<\infty \right\}.\end{equation}

As another consequence of Lemma \ref{growth} we get the following inversion
formula (see \cite{V-N} for the case $\alpha =0$).

\begin{corollary}
\label{inversion}Let $x\in X_{\alpha }$ and $h>s_{\alpha }(A)$.
If $F_{t}$ is defined as in Lemma \ref{inv}, then
\begin{equation*}
F_{t}(x)=\frac{1}{2\pi i}(C,1)\int_{\func{Re}\lambda =h}e^{\lambda
t}R(\lambda )xd\lambda
\end{equation*}
for all $t\in \mathbb{R}$.
\end{corollary}

\begin{proof}
If $h\geq \omega $, then our statement is the ordinary inversion formula
(see Lemma \ref{inv}). Otherwise, by the resolvent identity, we have
\begin{equation*}
R(u+iv)x=(1-(u-\omega )R(u+iv))A_{\omega }^{-\alpha }R(\omega +iv)A_{\omega
}^{\alpha }x\text{,}
\end{equation*}
for all $u$, $h\leq u\leq \omega $. So, in view of Lemma \ref{growth}, $%
\underset{v\rightarrow \infty }{\lim }R(u+iv)x=0$ uniformly in $u\in
[h,\omega ]$. Then, by the Cauchy Theorem, we get
\begin{eqnarray*}
\frac{1}{2\pi i}(C,1)\int_{\func{Re}\lambda =h}e^{\lambda t}R(\lambda
)xd\lambda &=&\frac{1}{2\pi i}(C,1)\int_{\func{Re}\lambda =\omega
}e^{\lambda t}R(\lambda )xd\lambda \\
&=&F_{t}(x)\text{.}
\end{eqnarray*}
\end{proof}

Let us recall the inequality $\omega _{\alpha }(\mathbf{T})\ge
s_{\alpha }(A)$ (see \cite{Weiss}). Now suppose $\omega >s_{\alpha }(A)$ and
$R(i\cdot +\omega )$ is a multiplier from $L_{p}(\mathbb{R},X_{\alpha })$ to $%
L_{p}(\mathbb{R},X)$. One can easily notice that the implication 1)$\Rightarrow
$4) of Theorem \ref{alpha} was proved individually for every $\rho $. Thus,
$$|\langle r_{\omega },\Phi \rangle |\leq |\langle \check{r}_{\omega },\hat{%
\Phi}\rangle |\leq K\Vert x\Vert _{\alpha }\Vert x^{*}\Vert \Vert \hat{\Phi}%
\Vert _{L_{1}}.$$ However, by Corollary \ref{inversion}, $\check{r}_{\omega
}(t)=e^{-\omega t}\langle x^{*},T_{t}x\rangle $, $t>0$, which implies $\Vert
T_{t}x\Vert \leq e^{\omega t}\Vert x\Vert _{\alpha }$. So, $\omega _{\alpha
}(\mathbf{T})\leq \omega $.

On the other hand, if $\omega >\omega _{\alpha }(\mathbf{T})$, then $\Vert
e^{-\omega t}T_{t}\Vert _{X_{\alpha }\rightarrow X}$ is exponentially
decaying. Consequently, the operator $M_{\omega }$, being a convolution with
the kernel $e^{-\omega t}T_{t}$, maps $L_{p}(\mathbb{R},X_{\alpha })$ into $%
L_{p}(\mathbb{R},X)$ as a bounded operator.
Thus, we have proved the following result.

\begin{theorem}
\label{constants}For any C$_{\text{0}}$-semigroup $\mathbf{T}$ on a Banach
space $X$, $\omega _{\alpha }(\mathbf{T})$ is the infimum over all $\omega
>s_{\alpha }(A)$ such that $R(i\cdot +\omega )$ is a multiplier from $L_{p}(%
\mathbb{R},X_{\alpha })$ to $L_{p}(\mathbb{R},X)$, for some $p$, $1\leq p<\infty $.
\end{theorem}

\begin{corollary}[\cite{Weiss,WW}]
If $X$ is a Hilbert space, then $\omega _{\alpha }(\mathbf{T})=s_{\alpha }(A)
$ for any strongly continuous semigroup $\mathbf{T}$ and $\alpha \geq 0$.
\end{corollary}
\noindent There are many results about properties of the constants $%
\omega _{\alpha }(\mathbf{T})$, $s_{\alpha }(A)$ and relations between them.
We refer the reader to paper \cite{WW} for a detailed exposition of the
subject.

We conclude this section by proving an $\alpha $-analogue of Perron's
Theorem, cf. \cite{LR}. Let us recall the classical result: a C$_{\text{0}}$-semigroup $%
\mathbf{T}$ with generator $A$ is hyperbolic if and only if for every $g\in
L_{p}(\mathbb{R},X)$ the following integral equation
\begin{equation}
u(\theta )=T_{\theta -\tau }u(\tau )+\int_{\tau }^{\theta }T_{\theta
-s}g(s)ds\text{, \ \ }\theta \geq \tau \text{,}  \label{mild}
\end{equation}
has unique solution in $L_{p}(\mathbb{R},X)$ (see, e.g. \cite[Theorem 4.33]{CL}%
).

In case of arbitrary $\alpha \geq 0$, we are looking for a necessary and
sufficient condition on $\mathbf{T}$, which provides existence and
uniqueness of solution to (\ref{mild}) in $L_{p}(\mathbb{R},X)$ for any given $%
g\in L_{p}(\mathbb{R},X_{\alpha })$. It turns out that the multiplier property
of $R(is)$ is the condition we need.

\begin{theorem}
\label{perron}Suppose $i\mathbb{R}\subset \rho (A)$. Then the following
assertions are equivalent:
\begin{itemize}
\item[1)] $R(i\cdot )$ is a multiplier from $%
L_{p}(\mathbb{R},X_{\alpha })$ to $L_{p}(\mathbb{R},X)$;
\item[2)] for every $g\in L_{p}(%
\mathbb{R},X_{\alpha })$ there exists a unique solution of (\ref{mild}) belonging
to $L_{p}(\mathbb{R},X)$.\end{itemize}
\end{theorem}

Before we prove the theorem let us state one auxiliary fact, see
 \cite{MiRaSc} or \cite[Prop.4.32]{CL}.
\begin{lemma}
\label{aux} A function $u$ is a solution of (\ref{mild}) if and only if $u\in D(\Gamma )
$ and $\Gamma u=-g$, where $\Gamma $ is the generator of the associated
evolution semigroup.
\end{lemma}

\begin{proof}
2) $\Rightarrow$ 1).
Denote by $L$ the linear operator that maps $g\in L_{p}(\mathbb{R},X_{\alpha })$ to the corresponding solution of (\ref{mild}). By the Closed
Graph Theorem, $L$ is bounded. We prove that actually $L=M_{0}$. Indeed, by
Lemma \ref{aux}, $Lg\in D(\Gamma )$ and $\Gamma Lg=-g$, for every $g\in
L_{p}(\mathbb{R},X_{\alpha })$. On the other hand, a straightforward
computation shows that if $g$ is a $C^{\infty }$-function with compact
support, then $M_{0}g\in D(\Gamma )$ and $\Gamma M_{0}g=-g$. Thus, $%
\Gamma (M_{0}g-Lg)=0$. However, if $\Gamma u=0$ for some $u\in D(\Gamma )$,
then again by Lemma \ref{aux}, $u$ is a solution of (\ref{mild})
corresponding to $g=0$. By the uniqueness, we get $u=0.$ So, $M_{0}g=Lg$
on a dense subspace of $L_{p}(\mathbb{R},X_{\alpha })$ and boundedness of $%
M_{0} $ is proved.

1) $\Rightarrow$ 2).  Suppose $M_{0}$ is bounded from $L_{p}(\mathbb{R},X_{\alpha })$ to $L_{p}(%
\mathbb{R},X)$. For a fixed $C^{\infty }$-function $g$ having compact support,
we show that $u=M_{0}g$ solves (\ref{mild}). Indeed, using \eqref{eq1}%
, we get
\begin{eqnarray*}
u(\theta )-T_{\theta -\tau }u(\tau ) &=&\int_{\mathbb{R}}R(is)\hat{g}%
(s)e^{is\theta }ds-\int_{\mathbb{R}}R(is)T_{\theta -\tau }\hat{g}(s)e^{is\tau
}ds \\
&=&\int_{0}^{\theta -\tau }T_{r}\int_{\mathbb{R}}e^{is(\theta -r)}\hat{g}(s)dsdr
\\
&=&\int_{0}^{\theta -\tau }T_{r}g(\theta -r)dr=\int_{\tau }^{\theta
}T_{\theta -r}g(r)dr\text{,}
\end{eqnarray*}
which is precisely (\ref{mild}).

Now suppose $g$ is an arbitrary function from $L_{p}(\mathbb{R},X_{\alpha })$.
Let us approximate $g$ by functions $(g_{n})$ of considered type. Then $%
u_{n}=M_{0}g_{n}$ converge to $u=M_{0}g$ in $L_{p}(\mathbb{R},X)$ and, without
loss of generality, pointwise on a set $E\subset \mathbb{R}$ with $mes\{\mathbb{R}%
\backslash E\}=0$. Thus, (\ref{mild}) is true for $u$, $g$ and all $\theta $
and $\tau $ from $E$. To get (\ref{mild}) for all $\theta $ and $\tau $, we
will modify $u$ on the set $\mathbb{R}\backslash E$. To this end, let us take a
decreasing sequence $(\tau _{n})\subset E$ such that $\lim \tau _{n}=-\infty
$. Observe that the functions $f_{n}(\theta )=T_{\theta -\tau _{n}}u(\tau
_{n})+\int_{\tau _{n}}^{\theta }T_{\theta -s}g(s)ds$ defined for $\theta
\geq \tau _{n}$ are continuous. Since $u=f_{n}=f_{m}$ on $(+\infty ,\max
(\tau _{n},\tau _{m})]\cap E$, we get $f_{n}=f_{m}$ everywhere in the
half-line $(+\infty ,\max (\tau _{n},\tau _{m})]$. Put $\tilde{u}$ to be $%
f_{n}$ on $(+\infty ,\tau _{n}]$. By the above, $\tilde{u}$ is a
well-defined function on all $\mathbb{R}$. Obviously, $u=\tilde{u}$ on $E$. Let
us show that $\tilde{u}$ satisfies (\ref{mild}). Indeed, for any $\theta
\geq \tau $ and $\tau >\tau _{n}$ we have
\begin{eqnarray*}
T_{\theta -\tau }\tilde{u}(\tau )+\int_{\tau }^{\theta }T_{\theta -s}g(s)ds
&=&T_{\theta -\tau }[T_{\tau -\tau _{n}}u(\tau _{n})+\int_{\tau _{n}}^{\tau
}T_{\tau -s}g(s)ds] \\
&&+\int_{\tau }^{\theta }T_{\theta -s}g(s)ds \\
&=&T_{\theta -\tau _{n}}u(\tau _{n})+\int_{\tau _{n}}^{\theta }T_{\tau
-s}g(s)ds=\tilde{u}(\theta ).
\end{eqnarray*}
\end{proof}

Clearly, assertion 1) in Theorem \ref{perron} is weaker than
condition 1) in Theorem \ref{alpha}. We do not know if they are equivalent.
In case $\alpha =0$, though, we can apply the resolvent identity to argue
that if $R(i\cdot )$ is a multiplier, then $R(i\cdot +\rho )$ is also a
multiplier for small values of $\rho $. So, by Theorem \ref{t5}, this
is equivalent to hyperbolicity of the semigroup $\mathbf{T}$, and our
statement turns into classical Perron's Theorem.

\section{An $\alpha $-analogue of hyperbolicity}

We begin with a discrete version of Theorem \ref{alpha} in the spirit of
\cite[Theorem 5]{LR}. Denote by $\func{Rg}T$ the range of an operator $T$.

\begin{theorem}
\label{discretealpha} Suppose $i\mathbb{Z}\subset \rho (A)$. Then the following
conditions are equivalent:

\begin{enumerate}
\item[1)]  $X_{\alpha }\subset \func{Rg}(I-T_{2\pi })$;

\item[2)]  The sum $(C,1)\sum_{k\in \mathbb{Z}}R(ik)x$ exists in $X$-norm for
all $x\in X_{\alpha }$;

\item[3)]  $\{R(ik)\}_{k\in \mathbb{Z}}$ is a multiplier from $L_{p}(\mathbb{T}%
,X_{\alpha })$ to $L_{p}(\mathbb{T},X)$ for some/all $1\leq p<\infty $;

\item[4)]  $\{R(ik)\}_{k\in \mathbb{Z}}$ is a multiplier from $L_{1}(\mathbb{T}%
,X_{\alpha })$ to $F(\mathbb{T},X)$, where $F$ is some quasi-normed function
lattice;

\item[5)]  There exists a constant $K>0$ such that
\begin{equation*}
|\langle r_{0},\Phi \rangle |=|\sum_{k\in \mathbb{Z}}r_{0}(k,x,x^{*})\Phi
(k)|\leq K\Vert x\Vert _{\alpha }\Vert x^{*}\Vert \Vert \check{\Phi}\Vert
_{L_{1}(\mathbb{T})}
\end{equation*}
holds for all $x\in X_{\alpha }$ , $x^{*}\in X^{*}$, and $\Phi \in C_{\infty
}(\mathbb{T})$.
\end{enumerate}
\end{theorem}

\begin{proof}
1)$\Leftrightarrow $2). Note that
\begin{equation*}
\frac{1}{2\pi }R(ik)(I-T_{2\pi })x=\frac{1}{2\pi }\int_{0}^{2\pi
}e^{-ikt}T_{t}xdt\text{,}
\end{equation*}
for all $x\in X$. So,
\begin{equation*}
\frac{1}{2\pi }(C,1)\sum_{k\in \mathbb{Z}}R(ik)(I-T_{2\pi })x=\frac{1}{2}%
(I+T_{2\pi })x.
\end{equation*}
Thus, 1) implies 2).

Now assume 1). Denote $S=\frac{1}{2\pi }(C,1)\sum_{k\in \mathbb{Z}}R(ik)x$.
Then
\begin{equation}
(\frac{1}{2}I+S)(I-T_{2\pi })x=(I-T_{2\pi })(\frac{1}{2}I+S)x=x\text{,}
\label{sleft}
\end{equation}
for all $x\in X_{\alpha }$, and 1) follows.

Following \cite{LR} we denote by $K$ the operator of convolution with the
semigroup, i.e.
\begin{equation*}
Kf(t)=\int_{0}^{2\pi }T_{s}f((t-s)[\func{mod}2\pi ])ds.
\end{equation*}
Clearly, $K$ is bounded on $L_{p}(\mathbb{T}, X)$ for all $1\leq p<\infty $ and $%
\alpha >0$. Now we define the discrete multiplier operator $L$ by the rule
\begin{equation*}
Lf(\theta )=\sum_{k\in \mathbb{Z}}R(ik)\hat{f}(k)e^{ik\theta },\quad
\theta\in[0,2\pi],
\end{equation*}
where $f$ is a trigonometric polynomial. One can check the identity
\begin{equation}
K=L(I-T_{2\pi })=(I-T_{2\pi })L.  \label{KLT}
\end{equation}

 By the assumption and the spectral mapping theorem for the point spectrum
 the operator $(I-T_{2\pi })$
is one-to-one. Suppose 1) holds. Then   $(I-T_{2\pi })$ has the left inverse
$U_1:=(I-T_{2\pi })_{left}^{-1}$ defined on $X_{\alpha }$. By the Closed Graph Theorem $U_1$
is bounded as an operator from $X_\alpha$ to $X$. Then (\ref{KLT}) says that $KU_{1}=L$
on trigonometric polynomials with values in $X_{\alpha }$. So, $L$ maps $%
L_{p}(\mathbb{T},X_{\alpha })$ into $L_{p}(\mathbb{T},X)$ for all $1\leq p<\infty $%
, which is what is stated in 3).

If the assertion in 3) is true only for some $p$, then as in the proof of Theorem \ref
{alpha}, $LA_{\omega }^{-\alpha }$ maps $L_{p}(\mathbb{T},X)$ into $L_{\infty }(%
\mathbb{T},X)$. By duality, $(LA_{\omega }^{-\alpha })^{*}$ maps $L_{1}(\mathbb{T}%
,X^{\odot })$ into $L_{q}(\mathbb{T},X^{\odot })$ and hence into $L_{\infty }(%
\mathbb{T},X^{\odot })$. So, $LA_{\omega }^{-\alpha }$ is a bounded operator
from $L_{1}(\mathbb{T},X)$ to $L_{\infty }(\mathbb{T},X)$, which proves 4) with $%
F=L_{\infty }(\mathbb{T},X)$.

Assume 4). Then for every $f\in C_{\infty }(\mathbb{T},X)$ there is a $\theta
\in [0,2\pi ]$ such that
\begin{equation*}
\left\| \sum_{k\in \mathbb{Z}}R(ik)\hat{f}(k)e^{ik\theta }\right\| \leq K\Vert
f\Vert _{1}\text{.}
\end{equation*}
Applying $T_{\theta }$ in the above inequality and using (\ref{eq1}) we
have:
\begin{equation*}
\left\| \sum_{k\in \mathbb{Z}}R(ik)\hat{f}(k)\right\| \leq K^{\prime }\Vert
f\Vert _{1}.
\end{equation*}
In particular, for $f=\check{\Phi}\otimes x$ the last inequality yields 5).

If 5) holds, then taking $\Phi =\sum_{n=1}^{N}\Phi _{n}\otimes x_{n}$, with $%
\Phi _{n}\in C_{\infty }(\mathbb{T})$ having disjoint supports we get the
following estimates
\begin{eqnarray*}
\Vert L\Phi \Vert _{1} &=&\frac{1}{2\pi }\int_{0}^{2\pi }\sup_{\Vert
x^{*}\Vert =1}\left| \sum_{n=1}^{N}\sum_{k\in \mathbb{Z}}r_{0}(k,x_{n},x^{*})%
\hat{\Phi}_{n}(k)e^{ik\theta }\right| d\theta \\
&\leq &\sum_{n=1}^{N}\frac{1}{2\pi }\int_{0}^{2\pi }\sup_{\Vert x^{*}\Vert
=1}|\langle r_{0},\widehat{\Phi (\cdot +\theta )}|d\theta \\
&\leq &\sum_{n=1}^{N}\Vert x_{n}\Vert _{\alpha }\Vert \Phi _{n}\Vert
_{1}=\Vert \Phi \Vert _{1}\text{.}
\end{eqnarray*}
So, we proved 3).

Finally, similarly to the proof of convergence of the $(C,1)$-integral in
Theorem \ref{alpha}, we can show that 2) follows from $L_{1}$-boundedness of
$LA_{\omega }^{-\alpha }$.
\end{proof}

In Theorem \ref{t5} we have proved that conditions
2) through 4), involving multipliers, are equivalent to the hyperbolicity of the semigroup, that
is, to a spectral property of ${\bf T}$. A natural question is to see if
the multipliers-type conditions 2) through 4) in Theorem \ref{alpha} are
equivalent to a spectral property  that could be formulated
in terms of ${\bf T}$ acting on the Banach space $X$, and {\em not}
in terms of a space of $X$-valued functions. Theorem \ref{discretealpha}
suggests that each of the conditions 2) through 4) in Theorem \ref{alpha}
implies that the inclusion $X_{\alpha }\subset \func{Rg}(zI-T_{2\pi })$
holds for all $z$ from some
annulus $\mathbb{A}$ containing $\mathbb{T}$. It turns out that this inclusion
alone is not equivalent to any of the conditions in Theorem
\ref{alpha}.
 Below we will find the needed complement, but let us first make
some observations.

Assume that $\mathcal{S}$ is some strip containing the imaginary axes and $%
\mathcal{S}\subset \rho (A)$. Suppose also that $X_{\alpha }\subset \func{Rg}%
(zI-T_{2\pi })$ for all $z$ from some annulus $\mathbb{A}$ containing $\mathbb{T}$%
. By the Point Spectrum Mapping Theorem $(zI-T_{2\pi })$ is one-to-one.
Thus, the left inverse operator $U_{z}=(zI-T_{2\pi })_{left}^{-1}:X_{\alpha
}\rightarrow X$ exists and is bounded by the Closed Graph Theorem. The
family $\mathbb{U}=\{U_{z}\}_{z\in \mathbb{A}}$ obeys the resolvent identity on
vectors from $X_{2\alpha }$. However, to prove analyticity, first of all one
needs uniform boundedness of $\mathbb{U}$. And that is the condition we are
looking for.

\begin{theorem}
\label{specalpha}Suppose there is a strip $\mathcal{S}$ such that $i\mathbb{R}%
\subset \mathcal{S}\subset \rho (A)$. Then any of the equivalent conditions
of Theorem \ref{alpha} holds if and only if there exists an annulus $\mathbb{A}$
containing $\mathbb{T}$ such that $X_{\alpha }\subset \func{Rg}(zI-T_{2\pi })$
for all $z\in \mathbb{A}$, and $\sup \{\Vert U_{z}\Vert _{X_{\alpha
}\rightarrow X}:z\in \mathbb{A}\}<\infty $.
\end{theorem}

\begin{proof}
Suppose that $M_{\rho }$ maps $L_{p}(\mathbb{R},X_{\alpha })$ into $L_{p}(\mathbb{R%
},X)$ for $|\rho |<2\rho _{0}$. By the Uniform Boundedness Principle $%
M_{\rho }$ are uniformly bounded for all $|\rho |<\rho _{0}$. Then, by
transference, $\{R(i(k+\xi )+\rho \}_{k\in \mathbb{Z}}$ is a multiplier
uniformly in $\xi \in [0,1)$ and $|\rho |<\rho _{0}$. In view of just proved
Theorem \ref{discretealpha} we get $X_{\alpha }\subset Rg(zI-T_{2\pi })$ for
all $z$ from some open annulus $\mathbb{A}$ containing $\mathbb{T}$. In order to
show uniform boundedness of $\mathbb{U}$, let us look at identity (\ref{sleft})
first. It shows, in particular, that $U_{1}=\frac{1}{2}+S$. Just like in the
second part of the proof of Theorem \ref{t5}, one can estimate $\Vert S\Vert
_{X_{\alpha }\rightarrow X}$ by the multiplier norm of $\{R(ik)\}$.
Rescaling gives the same conclusion for all $U_{z}.$ Since norms of the
corresponding multipliers are uniformly bounded, the desired result is
proved.

Now let us prove the converse statement.

Clearly, the family $\mathbb{U}=\{U_{z}\}_{z\in \mathbb{A}}$ obeys the resolvent
identity on vectors from $X_{2\alpha }$. Since, in addition, it is bounded,
the mapping $\ z\rightarrow U_{z}$ is strongly continuous on vectors from $%
X_{2\alpha }$ and, hence, on all $X_{\alpha }$. Again by the resolvent
identity, $U_{z}$ is strongly analytic on $X_{2\alpha }$. Since for any $%
x\in X_{\alpha }$, $U_{z}x$ is the uniform limit of a sequence $U_{z}x_{n}$
with $x_{n}\in X_{2\alpha }$, $U_{z}x$ is analytic.

It suffices to show that the integral $G(t)x=(C,1)\int_{\mathbb{R}%
}R(is)xe^{ist}ds$ converges for all $x\in X_{\alpha }$ and there exists a $%
\beta >0$ such that $\left\| G(t)x\right\| \leq Ke^{-\beta |t|}\left\|
x\right\| _{\alpha }$.

So, let us fix $x\in X_{\alpha }$ and $t\in \mathbb{R}$. Then for any $s\in
\mathbb{R}$, $x=(e^{2\pi is}-T_{2\pi })U_{e^{2\pi is}}x$. Thus
\begin{eqnarray*}
R(is)xe^{ist} &=&R(is)(I-e^{-2\pi is}T_{2\pi })U_{e^{2\pi
is}}xe^{ist}e^{2\pi is} \\
&=&e^{(2\pi +t)is}\int_{0}^{2\pi }e^{-ris}T_{r}U_{e^{2\pi is}}xdr.
\end{eqnarray*}
From this we get
\begin{eqnarray*}
G(t)x &=&\lim_{N\rightarrow \infty }\int_{0}^{2\pi
}T_{r}\int_{-N}^{N}U_{e^{2\pi is}}xe^{(2\pi +t-r)is}(1-\frac{|s|}{N})dsdr \\
&=&\lim_{N\rightarrow \infty }\int_{0}^{2\pi }T_{r}\int_{0}^{1}U_{e^{2\pi
is}}xe^{(2\pi +t-r)is}\sum_{n=-N}^{N}e^{(t-r)in}(1-\frac{|s+n|}{N})dsdr \\
&=&\lim_{N\rightarrow \infty }\int_{0}^{2\pi }T_{r}\int_{0}^{1}U_{e^{2\pi
is}}xe^{(2\pi +t-r)is}\sum_{n=-N}^{N}e^{(t-r)in}(1-\frac{|n|}{N})dsdr \\
&=&\lim_{N\rightarrow \infty }\int_{0}^{2\pi }T_{r}\int_{0}^{1}U_{e^{2\pi
is}}xe^{(2\pi +t-r)is}\mathcal{F}_{N}(-r+t\func{mod}2\pi )dsdr,
\end{eqnarray*}
where $\mathcal{F}_{N}$ is the Fej\'{e}r kernel. Passing to limit inside the
integral we get
\begin{eqnarray*}
G(t) &=&\frac{1}{2}\left[ T_{2\pi }\int_{0}^{1}U_{e^{2\pi
is}}xe^{ist}ds+\int_{0}^{1}U_{e^{2\pi is}}xe^{is(2\pi +t)}ds\right]  \\
&=&\frac{1}{2}\left[ T_{2\pi }\frac{1}{2\pi i}\int_{\mathbb{T}}z^{t}U_{z}xdz+%
\frac{1}{2\pi i}\int_{\mathbb{T}}z^{2\pi +t}U_{z}xdz\right] ,
\end{eqnarray*}
if $t=0\func{mod}2\pi $. And
\begin{eqnarray*}
G(t) &=&T_{t\func{mod}2\pi }\int_{0}^{1}U_{e^{2\pi is}}xe^{(2\pi +t-t\func{%
mod}2\pi )is}ds \\
&=&T_{t\func{mod}2\pi }\frac{1}{2\pi i}\int_{\mathbb{T}}z^{2\pi +t-t\func{mod}%
2\pi }U_{z}xdz,
\end{eqnarray*}
otherwise. In either case, replacing $\mathbb{T}$ by $(1+\varepsilon )\mathbb{T}$,
if $t<0$, or by $(1-\varepsilon )\mathbb{T}$, otherwise, we get the desired
exponential decay.
\end{proof}

\section{Strong $\alpha $-hyperbolicity\label{alphahypersec}}

In this section we introduce yet another notion of $\alpha $-hyperbolicity
for strongly continuous semigroups.
The spectral property we considered in the previous section, though strong
enough, fails to produce any splitting projection, which is so natural in
the case $\alpha >0$. Therefore, we investigate a notion of strong $\alpha $%
\textit{-}hyperbolicity, in which we force such a projection to exist.

\begin{definition}
\strut \label{hyperdef}A C$_{\text{0}}$-semigroup $\mathbf{T}=(T_{t})_{t\geq
0}$ is said to be strongly $\alpha $\textit{-hyperbolic} if there exists a
projection $P$ \ on $X$, called splitting, such that $PT_{t}=T_{t}P$, $t\geq
0$ and the following two conditions hold:

\begin{enumerate}
\item  $\omega _{\alpha }(\mathbf{T|}_{\func{Im}P})<0$;

\item  the restriction of $\mathbf{T}$ on $\func{Ker}P$ is a group, and $\omega _{\alpha }(%
\mathbf{T}^{-1}\mathbf{|}_{\func{Ker}P})<0$, where $\mathbf{T}%
^{-1}=(T_{-t}|_{\func{Ker}P})_{t\geq 0}$.
\end{enumerate}

The function $G(t)$ defined as in Definition \ref{Greenfn} is called \textit{%
the Green's function} corresponding to the $\alpha $-hyperbolic semigroup $%
\mathbf{T}$.
\end{definition}

It is an immediate consequence of the definition that Green's function
exponentially decays at infinity on vectors from $X_{\alpha }$.

Now we prove an analogue of Theorem \ref{t5} for $\alpha $-hyperbolic
semigroups.

\begin{theorem}
\strut \label{alphachar}A semigroup $\mathbf{T}$ is $\alpha $-hyperbolic if
and only if one of the equivalent conditions of Theorem \ref{alpha} is
satisfied and the operator $$G(t)x=\frac{1}{2\pi }(C,1)\int_{\mathbb{R}}R(is)xe^{ist}ds$$
 has a continuous extension to all of $X$ for each $t\in
\mathbb{R}$.

If this is the case, $G(t)$ represents the Green's function. Furthermore,
the splitting projection is unique and given by
\begin{equation}
P=\frac{1}{2}I+G(0).  \label{proj}
\end{equation}
\end{theorem}

\begin{proof}
Let us prove necessity.

If $\mathbf{T}$ is $\alpha $-hyperbolic, then there is a splitting
projection $P$. Suppose $x\in (\func{Im}P)_{\alpha }$. Then by Corollary \ref
{inversion} applied to the semigroup $\mathbf{T}|_{\func{Im}P}$, $\ $we have
$F_{t}(x)=G(t)x$. In particular, $Px=x=\frac{1}{2}x+G(0)x$. On the other
hand, if $x\in (\func{Ker}P)_{\alpha }$, then by the same reason, $\tilde{F}%
_{t}(x)=\frac{1}{2\pi }(C,1)\int_{\mathbb{R}}R(is,-A)xe^{ist}ds=-G(-t)x$, where
$\tilde{F}_{t}(x)$ is defined by
\begin{equation*}
\tilde{F}_{t}(x)=\begin{cases} T_{-t}x, & t>0\\ \frac{1}{2} x, & t=0\\ 0, &
t<0\end{cases}\text{.}
\end{equation*}

So, $Px=0=\frac{1}{2}x+G(0)x$. Since $X_{\alpha }=(\func{Im}P)_{\alpha }+(%
\func{Ker}P)_{\alpha }$ is dense in $X$, this shows that $G(t)$ continuously
extends to all of $X$ and equality (\ref{proj}) is true. The uniqueness of $%
P $ follows automatically from (\ref{proj}).

Since the $(X_{\alpha }\rightarrow X)$-norm of $G(t)$ is exponentially
decreasing and $M_{0}(\Phi )=G*\Phi $ \ for all $\Phi \in \mathcal{S}$, $%
M_{0}$ is bounded from $L_{1}(\mathbb{R},X_{\alpha })$ to $L_{1}(\mathbb{R},X)$.
To show boundedness of $M_{\rho }$, it is enough to notice that if $\mathbf{T%
}$ is $\alpha $-hyperbolic, then the scaled semigroup $e^{\rho \cdot }%
\mathbf{T}$ is also $\alpha $-hyperbolic, for small values of $\rho $.

Now we prove sufficiency.

Let us introduce the operator $P=\frac{1}{2}I+G(0)$. Since Theorem \ref{t5}
is valid, and hence formulas (\ref{Green}) in Remark \ref{greenrem} are
true, the norm of $T_{t}$ on $P(X_{\alpha })$ is exponentially decaying.
Consequently, by the ordinary inversion formula for Laplace transform, we
get $G(0)x=\frac{1}{2}x$, for all $x\in P(X_{\alpha })$. This implies $%
P^{2}=P$ on all $X$ , in view of the continuity of $P$. So, $P$ is a
projection.

Obviously, $PT_{t}=T_{t}P$. On the other hand, since $P(X_{\alpha })=(\func{%
Im}P)_{\alpha }$, we have $\omega _{\alpha }(\mathbf{T|}_{\func{Im}P})<0$ and
condition 1 of Definition \ref{hyperdef} is proved.

To show invertibility of $T_{t}$ on $\func{Im}(I-P)$, we apply formula (\ref
{Green}). It implies that $\Vert G(t)\Vert \Vert T_{-t}x\Vert \geq \Vert
x\Vert $, for $x$ in $\func{Im}(I-P)$, and hence, $T_{-t}|_{\func{Im}(I-P)}$
is invertible. Another application of (\ref{Green}) and the second part of
Theorem \ref{t5} proves condition 2 in Definition \ref{hyperdef}.
\end{proof}

\end{document}